\title{Universal Cycles of Discrete Functions}
\author{\small Britni LaBounty-Lay\\
\small East Tennessee State University\\
\and
\small Ashley Bechel \\
\small East Tennessee State University \\
\and
\small Anant P.~Godbole \\
\small East Tennessee State University \\
\small \tt{godbolea@etsu.edu} }
\documentclass [10 pt]{article}
\usepackage{amsfonts}
\usepackage[dvips]{graphicx}
\usepackage{amsthm}
\usepackage{amsmath}
\usepackage{amssymb}
\date{}
\begin{document}
\hsize=4.5truein
\vsize=7.2truein
\def\diam{{\rm diam}}
\def\ess{\rm ess}
\def\p{{\mathbb P}}
\def\ep{\varepsilon}
\def\P{{\rm Po}}
\def\cf{{\cal F}}
\def\cl{{\cal L}}
\def\e{{\mathbb E}}
\def\v{{\mathbb V}}
\def\l{\lambda}
\def\ll{{\ell_n}}
\def\a{{\alpha_n}}
\def\ph{\varphi(n)\sqrt n}
\def\dist{{\rm dist}}
\def\lr{\left(}
\def\rr{\right)}
\def\cd{\cdot}
\def\ts{\thinspace}
\def\lc{\left\{}
\def\rc{\right\}}
\def\qed{\vbox{\hrule\hbox{\vrule\kern3pt\vbox{\kern6pt}\kern3pt\vrule}\hrule}}
\newcommand{\hyp}{\mathcal{H}}
\newcommand{\hypp}{\mathcal{K}}
\newcommand{\Nu}{$\begin{Large}$\nu$\end{Large}$}
\newcommand{\ignore}[1]{}
\newtheorem{thm}{Theorem}
\newtheorem{result}{Result}
\newtheorem{lemma}[thm]{Lemma}
\newtheorem{cor}[thm]{Corollary}
\newtheorem{prop}[thm]{Proposition}
\maketitle
\begin{abstract}
A connected digraph in which the in-degree of any vertex equals
its out-degree is Eulerian; this baseline result is used as the basis of
existence proofs for universal cycles (also known as deBruijn cycles or  $U$-cycles) of several combinatorial objects.  
We present new results on the existence of universal cycles of certain
classes of functions.  These include onto functions, and 1-inequitable sequences on a binary alphabet.  In each case the connectedness of the underlying graph is
the non-trivial aspect to be established.
\end{abstract}
\section{Introduction and Terminology}

\medskip
\noindent
\bfseries \slshape (Informal) Definition 1
\normalfont
A universal cycle of a combinatorial object is a cyclic and ``efficient" listing of all the values of that object with no repetition.

\medskip
\noindent
\bfseries \slshape Example 1.
3-letter words on the binary alphabet.  

\normalfont \noindent The cycle
$$11100010$$
covers all the possible ``words", namely 111, 110, 100, 000, 001, 010, 101, and 011 in the smallest sized sequence, namely of length 8.

\medskip
\noindent
\bfseries \slshape Example 2
All 3-subsets of $\{1,2,3,4,5,6,7,8\}$.

\normalfont\noindent Note that the U-cycle below is obtained (\cite{hurlbert}) by starting with the 7 numbers 1356725 and then adding 5 mod 8 successively to get the sequence 
$$1356725\enspace 6823472\enspace 3578147\enspace 8245614\enspace 5712361\enspace 2467836\enspace 7134582\enspace 4681258.$$
This U-Cycle (constructed by Glenn Hurlbert) of length ${8\choose 3}=56$ is read as follows: $\{1,3,5\},\{3,5,6\},\{5,6,7\}\ldots\{2,5,8\},\{5,8,1\},\{8,1,3\}.$  In fact Hurlbert shows in \cite{hurlbert} that for $k=3,4,6$, there is an integer $n_0(k)$ such that for $n\ge n_0(k)\Rightarrow$ U-cycles of $k$-subsets of $[n]=\{1,2,\ldots,n\}$ exist.

The following result is found as standard fare in most graph theory texts, see e.g. Theorem 1.4.24 in \cite{west}
\begin{thm}
A connected digraph is Eulerian if and only if the in-degree of each vertex equals its out degree.
\end{thm}
The next key result in the area of U-cycles is known as DeBruijn's theorem; see e.g. Theorem 1.4.26 in \cite{west} for a special case.
\begin{thm}
U-Cycles of $k$-letter words on an $n$-letter alphabet exist for all $k$ and $n$.
\end{thm}

\medskip

\noindent{\bf Proof.}  Identical to the proof of Theorem 1.4.26 in \cite{west}.  The key idea is to create a graph $G$ with vertex set that consists of all $k-1$ letter words on the $n$ letter alphabet, i.e. of length one less than the original word length, with a directed  edge being drawn from vertex $v_1$ to vertex $v_2$ if the last $k-2$ letters of the ``word" $v_1$ coincide with the first $k-2$ letters of $v_2$.  Here and in {\it all} the other situations we consider in this paper, we {label} the edge with the concatenated $k$-letter word thus formed, e.g. for $n=26,k=4$ the edge from $CAT$ to $ATE$ is labeled as $CATE$.  It is clear that each vertex has in-degree and out-degree equal to $n$.  Theorem 1 tells us that $G$ is Eulerian; if the Eulerian circuit happens to be, e.g., $CATE, ATEZ, TEZO,\ldots,RCAT$ then the corresponding U-cycle is $CATEZO\ldots R$.

The following result was proved in \cite{brad}
\begin{thm}
A U-Cycle of $1-1$ functions from $\left\{1,...,k\right\} \rightarrow \left\{1,...,n\right\}$ exists if and only if $n>k$; these are merely permutations of $n$ objects taken $k$ at a time, or, alternatively, $k$-letter words on $[n]$ in which no letter repeats.
\end{thm}
When $k=n$ it turns out that the underlying graph is not connected and thus a U-cycle cannot exist.  For example for $n=k=3$ the graph with vertex set consisting of two-letter words on $\{1,2,3\}$ decomposes into the cycles
\[12\rightarrow23\rightarrow31\]
and
\[21\rightarrow13\rightarrow32.\]
Notice that in Theorem 3, the induced edge labels consist of the objects we seek to build a U-cycle of, namely one-to-one functions, whereas the vertices are {\it also} one to one functions with domain (word length) of size $k-1$.  This will be in marked contrast to our Theorems 4 and 5, in which the vertices will often represent different kinds of combinatorial objects.  The proof of Theorem 3 utilizes Theorem 1 but is non constructive.  Knuth \cite{knuth} raised the question of when a U-cycle of one to one functions can be explicitly constructed and the first such effort appears to be, for $k=n-1$, due to Ruskey and Williams \cite{ruskey}.

In \cite{cdg}, Chung, Diaconis, and Graham {\it do} consider the case of $k=n$ and deal with U-cycles of all the $n!$ permutations on $[n]$ but obviously in a different sense than what we consider above.  Now, one to one functions constitute a restricted class of all functions from $[k]$ to $[n]$.  In this paper, we focus on this line of inquiry, turning our attention first from one to one to {\it onto} functions and then, when $n=2$, to {\it 1-inequitable} functions, which we define below.
\section{Results}
\noindent
\bfseries \slshape Definition 2
\normalfont  A function $f:[k]\rightarrow[n]$ is said to be {\it almost onto} if $\vert[n]\setminus{\rm Range}(f)\vert=1$

\medskip
\noindent
\bfseries \slshape Definition 3
\normalfont
A binary word of even length $k\ge4$ is said to be {\it equitable} if it consists of $k/2$ zeros and $
k/2$ ones.

\medskip
\noindent
\bfseries \slshape Definition 4
\normalfont
A binary word of odd length $k\ge 3$ is said to be {\it 1-inequitable} if it consists of $ \lfloor k/2\rfloor$ zeros and $
\lceil k/2\rceil$ ones -- or vice versa.

\medskip
\noindent
\bfseries \slshape Definition 5
\normalfont
A binary word of even length $k$ is said to be {\it 2-inequitable} if the numbers of ones and zeros differ by two.


\begin{thm}

A U-Cycle of $onto$ functions from $\left\{1,...,k\right\} \rightarrow \left\{1,...,n\right\}$ exists iff $k>n$

\end{thm}

\medskip\noindent{\bf Proof.}  We have already seen that a U-cycle of onto functions cannot exist when $k=n$; assume, henceforth that $k>n$.  The vertices of the underlying graph consist of {\it certain kinds} of functions from $\left\{1,...,k-1\right\} \rightarrow \left\{1,...,n\right\}$ -- these are functions such that the corresponding concatenated edge label consists of an onto function from $\left\{1,...,k\right\} \rightarrow \left\{1,...,n\right\}$.  A moment's reflection reveals that the only allowable such vertices are those corresponding to {\it onto} and {\it almost onto} functions on $\{1,2,\ldots,k-1\}$.  Moreover, there is a dichotomy between the degree structure of these two classes of vertices:  If $v$ is onto, its indegree and outdegree both equal $n$, however $i(v)=o(v)=1$ for almost onto vertices, where $i(v)$ and $o(v)$ represent the indegree and outdegree of $v$ respectively.  What is critical, though, is that $i(v)=o(v)$ for each $v$.  To invoke Theorem 1, we need to show that our graph is connected, i.e., that there is a path from $u$ to $v$, no matter what {\it kind} of vertices $u,v$ happen to be.  We start by assuming that both $u$ and $v$ are onto functions from $\left\{1,...,k-1\right\} \rightarrow \left\{1,...,n\right\}$ and will now exhibit the fact that there is a path between them.  The idea is simple.  Set $k-1=M$.  Suppose that we seek to build a path between $a=a_1a_2\ldots a_M$ and $b=b_1b_2\ldots b_M$, where $a$ and $b$ are both onto.  We start by ``building" the sequence $b$ to the extent that it is ``legal,", i.e. by traversing the trail
\begin{equation}a_1a_2\ldots a_M\rightarrow a_2\ldots a_Mb_1\rightarrow\ldots\rightarrow a_r\ldots a_Mb_1\ldots b_{r-1}\end{equation}to the extent possible, i.e. for some $r\ge2$. At this point the word we have arrived at must be almost onto, for if it were onto we could continue the process.  Now we claim that either two of the $a$s in the last sequence in (1) must be the same, or else one of the $a$s must equal one of the $b$s.  For, if not, the $a$s must all be distinct, and, since the set of $a$s and the set of $b$s are disjoint it would be impossible for $b$ to be onto.  Let $a_{s-1}$ be the first $a$, from the left, that is a ``duplicate".  

Now the outdegree of the last word in (1) is one; we continue to add the only allowable letter as we build the trail through a series of almost onto words, {\it ultimately arriving at an onto word $c$ as follows:}
\[a_r\ldots a_Mb_1\ldots b_{r-1}\rightarrow a_{r+1}\ldots a_M b_1\ldots b_{r-1}\diamondsuit_1\ldots\]\[\rightarrow c=a_s\ldots a_Mb_1b_2\ldots b_{r-1}\diamondsuit_1\ldots\diamondsuit_{s-r},\]
where the $\diamondsuit_j$'s are those letters forced to be added on by the fact that we are building the trail through almost onto words, thus having only one choice for the next vertex in the trail.  But $c$ is onto, so that we may next travel from it to a cyclic version
\[c^*=\diamondsuit_1\ldots\diamondsuit_{s-r}a_s\ldots a_Mb_1\ldots b_{r-1}\]
which is also onto.  But now we may travel to $\diamondsuit_2\ldots\diamondsuit_{s-r}a_s\ldots a_Mb_1\ldots b_{r}$ and we are thus able to build one more letter, namely $b_r$, in our quest to travel to the sequence $b$ as begun in (1).  We iteratively continue this process until the word $b$ is reached.  

Let us illustrate the above process by an example:  Suppose $M=n=5$ and we wish to exhibit a path from $a=13425$ to $b=41235$.  We first add on the first two letters of $b$ as follows:
\[13425\rightarrow34254\rightarrow42541.\]
Next we travel from 42541 to an onto word as follows
\[42541\rightarrow25413,\]
which we cycle around until the ``41" segment is at the tail as follows
\[25413\rightarrow54132\rightarrow41325\rightarrow13254\rightarrow32541\]
which allows us to travel to 25412, whence $b$ may be reached easily:
\[25412\rightarrow54123\rightarrow41235.\]

Continuing with the proof of Theorem 4, there are three more cases that need to be considered to establish that $G$ is connected, namely that there is a path between an almost onto $a$ and an onto $b$ (or vice versa), or between two almost onto vertices $a,b$.  Let $a$ be almost onto and $b$ onto.  Let $a=a_1a_2\ldots a_M$ and let $a_s$ be the first letter that is represented twice in $a$.  We then proceed from $a$ as follows:
\[a_1a_2\ldots a_M\rightarrow\ldots\rightarrow a_{s-1}\ldots a_M\diamondsuit_1\ldots\diamondsuit_{s-2},\]
from which we travel to an onto word $c$ in a single step by reintroducing the missing letter.  Finally we can find a path from $c$ to $b$ as in the first part of the proof.  If $a$ is onto and $b$ is almost onto, then the strategy is inverse to the one indicated above; we first traverse a path from $a$ to a ``logical" onto vertex $c$, from which the path to $b$ is easy to establish.  Specifically, if $b_s$ is the first letter, from the right of the word $b$, that is represented twice, then we ``backtrack" from $b$ to $c$ as follows:
\[b=b_1b_2\ldots b_M\leftarrow\clubsuit_1b_1\ldots b_{M-1}\leftarrow\ldots\leftarrow\clubsuit_{M-s}\ldots\clubsuit_{1}b_1\ldots b_s\]
\[\leftarrow c=\heartsuit\clubsuit_{M-s}\ldots\clubsuit_{1}b_1\ldots b_{s-1},\]
where $c$, by construction, is onto.  A path from $a$ to $c$ is found as in Case 1 of the theorem.  Finally, the fourth case is proved by combining Cases 2 and 3.  This completes the proof.

\bigskip

Our next result deals with a special class of onto functions.  If the alphabet is binary, onto functions consist of all binary sequences except for $(1,1,\ldots,1)$ and $(0,0,\ldots,0)$, which makes the situation rather uninteresting since it is very close to DeBruijn's theorem.  We thus focus, in the binary case, on a {\it smaller} class of onto functions that we show admits a U-cycle.

\begin{thm}

A U-Cycle of $1-inequitable$ functions from $\left\{1,...,k\right\} (k\equiv1\mod2)\rightarrow \left\{1,0\right\}$ exists, while U-cycles of equitable functions from $\left\{1,...,k\right\}$ 

\noindent $(k\equiv0\mod2)\rightarrow \left\{1,0\right\}$ do not exist.
\end{thm}

\medskip\noindent
{\bf Proof.}  
Let us first prove that U-cycles of equitable binary functions do not exist.  For small values of $k$, say $k=4$, the fact that the underlying graph is disconnected is easy to see.  Vertices consist of binary words of length $k-1$ that are 1-inequitable and $i(v)=o(v)=1$ for each vertex $v$.  The graph decomposes for $k=4$ into the 2 cycles 
\[C_1=110\rightarrow100\rightarrow001\rightarrow011\rightarrow110\]
and \[C_2=101\leftrightarrow010.\]  In general, we might ask how many cycles $a_k$ the graph $G$ with vertex set consisting of 1-inequitable binary functions decomposes into; the above shows that $a_4=2$ and there is one cycle of length 4, written, in terms of edges as
\[1100\rightarrow1001\rightarrow0011\rightarrow0110\rightarrow1100\] and another of length 2, namely
\[1010\rightarrow0101\rightarrow1010.\]  The solution to the above question is rooted in the number of divisors of $k$.  If a word is $p$-periodic, then it will generate a cycle of length $p$.  Let $k_1, k_2,\ldots k_r$ denote the even divisors of the even integer $k$ and let $b_r$ denote the number of cycles with length $r$.  We clearly have 
\[a_k=b_k+b_{k_1}+\ldots b_{k_r},\eqno(2)\]
and 
\[kb_k+\sum_{j=1}^r k_jb_{k_j}={{k}\choose{k/2}},\]
or, in more useful terms,
\[b_k=\frac{{{k}\choose{k/2}}-\sum_{j=1}^rk_jb_{k_j}}{k}.\eqno(3)\]
Equations (2) and (3) lead to the sequences
\[\{a_{2k}\}=1,2,4,10,26,80,246,810,2704,9252,32066,112720\ldots\]
and
\[\{b_{2k}\}=1,1,3,8,25,75,245,800,2700,9225, 32065,112632,\ldots\]
which may be found in Neil Sloan's website of integer sequences 

\centerline{http://www.research.att.com/$\sim$njas/sequences} 
\noindent as Sequences A003239 and A022553 respectively, indicating that (albeit in a slightly different context) this problem had been previously solved.  

We thus move to the main part of the proof, namely showing that a U-cycle of 1-inequitable functions exists.  In this case, the underlying graph has vertex set that consists either of equitable sequences of length $k-1$ ($i(v)=o(v)=2$) or 2-inequitable sequences of length $k-1$ ($i(v)=o(v)=1$).  We next establish connectedness.  Assume first that we wish to traverse a path from $a=a_1a_2\ldots a_M$ to $b=b_1b_2\ldots b_M$, both equitable sequences.  We start by building $b$ as in the proof of Theorem 4 as follows:
\[a_1a_2\ldots a_M\rightarrow\ldots\rightarrow a_r\ldots a_Mb_1\ldots b_{r-1},\eqno(4)\]
where the last word in the above chain is 2-inequitable and has, without loss of generality, $\lfloor k/2\rfloor$ zeros and $\lceil k/2\rceil$ ones.  We thus add a ``0" to the chain to reach the vertex $a_{r+1}\ldots a_Mb_1\ldots b_{r-1}0$, which may be either equitable or 2-inequitable.  Assuming it is the latter,  we  travel from here through a possibly empty set of 2-inequitable vertices until ultimately we reach an equitable vertex (this will occur in the step immediately after the first ``1" among the $a$s above is reached.  Note that one of the $a$s {\it must} be a 1 since the word $b$ has been assumed to be equitable and thus has $M/2$ ones; the word in (4), however, has $(M/2)+1$ ones, one of which must come from the $a$ segment.)

Let the equitable vertex thus reached be denoted by 
$$\Delta=c_{\ell+1}\ldots c_M b_1b_2\ldots b_{r-1}c_r\ldots c_\ell.$$  We next travel through cyclic versions of $\Delta$, ending at $c_r\ldots c_Mb_1\ldots b_{r-1}$, which permits us to travel to $c_{r+1}\ldots c_M b_1\ldots b_{r}$, one step beyond what we had achieved in (4).  This algorithm is implemented until the word $b$ is reached.  The proof of the other three cases, namely establishing a path between (i) an equitable and a 2-inequitable vertex; (ii) a 2-inequitable vertex and an equitable vertex; and (iii) two 2-inequitable vertices is similar to that in Theorem 4 and is omitted.

\section{Open Problems}  The kinds of open problems that this paper raises concern, for example, the existence of universal cycles of certain kinds of inequitable functions when we are no longer restricted to a binary alphabet; or the existence of U-cycles of functions with growth-like conditions (possibly discrete Lipschitz-type conditions).  There are many possibilities.
\section{Acknowledgments}
This work 
forms part of the undergraduate research project of Ashley Bechel and Britni LaBounty-Lay conducted under the supervision of Anant Godbole, who was supported by NSF Grant DMS-0552730.


\begin{thebibliography}{99}

\bibitem{cdg} F.~Chung, P.~Diaconis, and R.~Graham (1992), ``Universal cycles for combinatorial structures," {\it Discrete Math.} {\bf 110}, 43--59.
\bibitem{hurlbert} Hurlbert, G. (1994), ``On universal cycles for $k$-subsets of an $n$-element set," {\it SIAM J. Discrete Math.} {\bf 7}, 598--604.
\bibitem{brad} Jackson, B. (1993), ``Universal cycles of $k$-subsets and $k$-permutations," {\it Discrete Math.} {\bf 117}, 114--150.
\bibitem{knuth} Knuth, D. (2005), {\it The Art of Computer Programming, Volume 4, Fascicle 2}, Pearson, NJ.
\bibitem{ruskey} Ruskey, F. and Williams, A. (2008), ``An explicit universal cycle for the $n-1$-permutations of an $n$-set," Talk at Napier Workshop.
\bibitem{west} West, D. (1996), {\it Introduction to Graph Theory}, Prentice Hall, New Jersey.
\end{thebibliography}
\end{document}